\documentclass[12pt]{amsart}
\usepackage{amscd,amssymb}
\usepackage[arrow,matrix]{xy}
\usepackage[plainpages,backref,urlcolor=blue]{hyperref}

\theoremstyle{plain}
\newtheorem{thm}[subsection]{Theorem}
\newtheorem{lem}[subsection]{Lemma}
\newtheorem{prop}[subsection]{Proposition}
\newtheorem{cor}[subsection]{Corollary}

\theoremstyle{definition}
\newtheorem{rk}[subsection]{Remark}
\newtheorem{definition}[subsection]{Definition}
\newtheorem{ex}[subsection]{Example}
\newtheorem{question}[subsection]{Question}

\numberwithin{equation}{section}
\newcommand{\OO}{{\mathcal O}}

\newcommand{\F}{{\mathcal F}}

\newcommand{\LL}{{\mathcal L}}

\newcommand{\al}{{\alpha}}
\newcommand{\be}{{\beta}}
\newcommand{\NN}{{\mathcal N}}

\newcommand{\Q}{\mathbb{Q}}

\newcommand{\C}{\mathbb{C}}

\newcommand{\PP}{\mathbb{P}}
\newcommand{\HH}{\mathbb{H}}

\DeclareMathOperator{\im}{im}
\DeclareMathOperator{\coker}{coker}

\DeclareMathOperator{\defect}{def}
\DeclareMathOperator{\codim}{codim}

\begin{document}
%\date{June 4, 2009}

\title [Koszul complexes and pole order filtrations]
{Koszul complexes and pole order filtrations}

\author[Alexandru Dimca]{Alexandru Dimca$^1$}
\address{Univ. Nice Sophia Antipolis, CNRS,  LJAD, UMR 7351, 06100 Nice, France. }
\email{dimca@unice.fr}

\author[Gabriel Sticlaru]{Gabriel Sticlaru}
\address{Faculty of Mathematics and Informatics,
Ovidius University,
Bd. Mamaia 124, 900527 Constanta,
Romania}
\email{gabrielsticlaru@yahoo.com }
\thanks{$^1$ Supported by Institut Universitaire de France}



\begin{abstract} We study the interplay between the cohomology of the Koszul complex of the partial derivatives of a homogeneous polynomial $f$ and the pole order filtration $P$ on the cohomology of the open set $U=\PP^n \setminus D$, with $D$ the hypersurface defined by $f=0$.
The relation is expressed by some spectral sequences, which may be used on one hand to determine the 
filtration $P$ in many cases for curves and surfaces, and on the other hand to obtain information about the syzygies involving the partial derivatives of the polynomial $f$. The case of a nodal hypersurface $D$ is treated in terms of the defects of linear systems of hypersurfaces of various degrees passing through the nodes of $D$. When $D$ is a nodal surface in $\PP^3$, we show that $F^2H^3(U) \ne P^2H^3(U)$ as soon as the degree of $D$ is at least $4$.
\end{abstract}

\maketitle

%\tableofcontents

\section{Introduction} \label{sec:intro}

Let $S=\C[x_0,...,x_n]$ be the graded ring of polynomials in $x_0,,...,x_n$ with complex coefficients and denote by $S_r$ the vector space of homogeneous polynomials in $S$ of degree $r$. 
For any polynomial $f \in S_N$ we define the {\it Jacobian ideal} $J_f \subset S$ as the ideal spanned by the partial derivatives $f_0,...,f_n$ of $f$ with respect to $x_0,...,x_n$. For $n=2$ we use $x,y,z$ instead of
$x_0,x_1,x_2$ and
$f_x,f_y,f_z$  instead of $f_0,f_1,f_2$, in the same way as in Eisenbud's book \cite{Eis}. 

We define the corresponding graded {\it Milnor} (or {\it Jacobian}) {\it algebra} by
\begin{equation} 
\label{eq1}
M(f)=S/J_f.
\end{equation}
The study of such Milnor algebras is related to the singularities of the corresponding projective hypersurface $D:f=0$, see \cite{CD}, as well as to the mixed Hodge theory of 
the hypersurface $D$ and of its complement $U=\PP^n \setminus D$, see the foundational article by Griffiths \cite{Gr} and also \cite{Dc}, \cite{DS2}, \cite{DSW}. For mixed Hodge theory we refer to \cite{PS}.

In fact, such a Milnor algebra can be seen (up to a twist in grading) as the first (or the last) homology (or cohomology) of the Koszul complex of the partial derivatives $f_0,...,f_n$ in $S$, see  \cite{CD} or  \cite{D1}, Chapter 6. As such, it is related to certain natural $E_1$-spectral sequences associated to the pole order filtration and converging to the cohomology of the complement $U$ introduced in \cite{Dc} and discussed in detail in \cite{D1}, Chapter 6.

In the second section we recall and improve the construction of these spectral sequences and show that 
they degenerate at the $E_2$-terms when all the singularities of $D$ are weighted homogeneous and $\dim D=1$, case when we use the notation $C$ for the curve $D$,
see Theorem \ref{thm2} $(iii)$. This result gives a positive answer in the curve case to an old conjecture by the first author, see the claim just before Remark (3.11) in \cite{Dc}. Such a degeneracy at the $E_2$-terms is shown to occur also for nodal surfaces, see Theorem \ref{thm7} $(i)$.

In the third section we assume $n=2$ and use this approach to determine the pole order filtration $P^*$ on the cohomology group $H^2(U)$ in a number of cases,
see Examples \ref{FP1}, \ref{FP2} and \ref{FP3}, the latter  being a new example where $F^2 \ne P^2$ on $H^2(U)$.
We also describe in Example \ref{FP0} these spectral sequences completely when $C$ is a nodal curve. 
 
In the forth section, we discuss the syzygies of nodal hypersurfaces. For instance we show that for a nodal curve there are no nontrivial relations
 \begin{equation} 
\label{rel1}
R_m:~~ af_x+bf_y+cf_z=0
\end{equation}
with $a,b,c$ homogeneous of degree $m<N-2$ and we describe completely the relations of degree $m=N-2$ in terms of the irreducible factors $f_j$ of $f$, see Theorem \ref{thm5}. (Note that $f_j$ has a different meaning for $n=2$ and for $n>2$). The vanishing part in Theorem \ref{thm5} was extended to nodal hypersurfaces of arbitrary dimension in  \cite{DSt3} , using a different approach.

\begin{definition}
\label{def}

For a hypersurface $D:f=0$ with isolated singularities we introduce three integers, as follows:

\noindent (i) the {\it coincidence threshold} $ct(D)$ defined as
$$ct(D)=\max \{q~~:~~\dim M(f)_k=\dim M(f_s)_k \text{ for all } k \leq q\},$$
with $f_s$  a homogeneous polynomial in $S$ of degree $N$ such that $D_s:f_s=0$ is a smooth hypersurface in $\PP^n$.

\noindent (ii) the {\it stability threshold} $st(D)$ defined as
$$st(D)=\min \{q~~:~~\dim M(f)_k=\tau(D) \text{ for all } k \geq q\}$$
where $\tau(D)$ is the total Tjurina number of $D$, i.e. the sum of all the Tjurina numbers of the singularities of $D$.

\noindent (iii) the {\it minimal degree of a nontrivial syzygy} $mdr(D)$ defined as
$$mdr(D)=\min \{q~~:~~ H^n(K^*(f))_{q+n}\ne 0\}$$
where $K^*(f)$ is the Koszul complex of $f_0,...,f_n$ with the grading defined in the next section.

\end{definition}

If a relation as in \eqref{rel1} is of minimal degree among the relations modulo the trivial relations \eqref{tij}, then one has
$mdr(D)=m$, i.e. our notion is the natural one. Moreover
it follows from \eqref{eq22} that one has 
\begin{equation} 
\label{REL}
ct(D)=mdr(D)+N-2.
\end{equation} 

By definition, it follows that for any such hypersurface $D$ which is not smooth, we have $N-2 \leq ct(D) \leq (n+1)(N-2)$ and using \cite{CD} we get $st(D) \leq (n+1)(N-2)+1$.
With these handy notation, we can state the following result, which is a consequence of the vanishings in Theorem \ref{thm5} obtained via Hodge theory, using the equation \eqref{eq22}.

\begin{thm}
\label{thm3}
Let $C:f=0$  be a nodal curve of degree $N$ in $\PP^2$.
Then one has $ct(C)\geq 2N-4$.
\end{thm}

Recall that Hilbert-Poincar\'e series of a graded $S$-module $E$ of finite type is defined by 
\begin{equation} 
\label{eq21}
HP(E)(t)= \sum_{k\geq 0} (\dim E_k)t^k
\end{equation}
and that we have
\begin{equation} 
\label{eq31}
HP(M(f_s))=  \frac{(1-t^{N-1})^{n+1}}{(1-t)^{n+1}}.
\end{equation}
In particular, if we set $T=T(n,N)=(n+1)(N-2)$, it follows that $M(f_s)_j=0$ for $j>T$ and $\dim M(f_s)_j=\dim M(f_s)_{T-j}$ for $0 \leq j \leq T$.

\medskip

Theorem \ref{thm3} determines the dimensions of $M(f)_q$ for all $q<2N-3$ in the case of a nodal curve $C$. The next dimension for such a curve is given by
\begin{equation} 
\label{dim2*}
\dim M(f)_{2N-3}=n(C)+\sum_{j=1,r}g_j=g+r-1.
\end{equation}
where $n(C)=\tau(C)$ is the total number of nodes of $C$ and $g_j$ are the genera of the normalizations of the irreducible components $C_j$ of $C$,
whose number is $r$ and 
\begin{equation} 
\label{dim1}
g=\frac{(N-1)(N-2)}{2}
\end{equation}
is the genus of the smooth curve $C_s:f_s=0$, see \eqref{dim2} and \eqref{dim3}. For more general curves we have the following relation between the Milnor algebra $M(f)$ and the geometry of $U$, consequence
of Theorem \ref{thm2}.

\begin{cor}
\label{corB1}
Let $C:f=0$ be a curve in $\PP^2$ of degree $N$, having only isolated weighted homogeneous singularities.
Then
$$\dim M(f)_{2N-3}+\dim P^2H^2(U)=2g+r-1=\dim H^2(U)+\tau (C)$$
where $g$ is as in \eqref{dim1}, $r$ is the number of irreducible components of $C$ and $\tau(C)$ is the total Tjurina number of $C$.

\end{cor}

\bigskip
 
For a highly singular curve $C$, we can have
much lower values for $ct(C)$ than that given by Theorem \ref{thm3}, namely one can have $ct(C)=N-2$ or $ct(C)=N-1$,
see Example \ref{lowdegrel}.

On the other hand, it follows from \eqref{dim2*} that for a nodal curve $C$ one has $ct(C)=2N-4$ if and only if $C$ is not irreducible (i.e. $r>1$).
One of the main results in \cite{DSt}, restated as the first equality in  
\eqref{dim2*}, implies that for a {\it rational nodal curve} (i.e. $g_j=0$ for $j=1,...,r$) one has $st(C)\leq 2N-3$. This yields the following.

\begin{cor}
\label{corA1}
For a rational nodal curve $C$, the Hilbert-Poincar\'e series $HP(M(f))$ is completely determined by the degree $N$ and the number of nodes $n(C)$.
In particular, $st(C)=2N-3$ unless $C$ is a generic line arrangement and then $st(C)=2N-4$ for $N>3$ and $st(C)=1$ for $N=3$.
\end{cor}
For the case of hyperplane arrangements, an interesting approach to the study of the Jacobian ideal $J_f$ is given in the recent paper \cite{DSSWW}.

\medskip

At the other extreme, there are nodal curves with $ct(C)=3N-6$, as implied by the description of the Hilbert-Poincar\'e series $HP(M(f))$ for any hypersurface having exactly one node given in Example \ref{exonenode}, $(i)$. 

To state the next result, we recall some notation. For a finite set of points $\NN \subset \PP^n$
we denote by 
$$\defect S_m(\NN)=|\NN| - \codim \{h \in S_m~~|~~ h(a)=0 \text{ for any } a \in \NN\},$$ 
the {\it defect (or superabundance) of the linear system of polynomials in $S_m$ vanishing at the points in $\NN$}, see \cite{D1}, p. 207. This positive integer is called the {\it failure of $\NN$ to impose independent conditions on homogeneous polynomials of degree $m$} in \cite{EGH}. In the fourth section we prove the following.

\begin{thm}
\label{onenode}
Let $D:f=0$ be a degree $N$ nodal hypersurface in $\PP^n$ and let $\NN$ denote the set of its nodes. Then
$$\dim H^{n}(K^*(f))_{nN-n-1-k}= \defect S_k(\NN) $$
for $0\leq k \leq nN-2n-1$  and $\dim H^{n}(K^*(f))_j=\tau(D)=|\NN|$ for $j\geq n(N-1)$.

In other words, 
$$\dim M(f)_{T-k}=\dim M(f_s)_k+\defect S_k(\NN) $$
for $0\leq k \leq nN-2n-1$, where $T=T(n,N)=(n+1)(N-2)$. In particular, $\dim M(f)_T=\tau(D)$, i.e. $st(D) \leq T$.
\end{thm}
Note that this Theorem determines the dimensions $\dim M(f)_j$ in terms of defects of linear systems for any $j\geq N-1$, i.e. for all $j$ since
the dimensions $\dim M(f)_j=\dim S_j$ for $j<N-1$ are well known. The last equality, namely $\dim M(f)_T=\tau(D)$, improves by one Corollary 9. in \cite{CD} in the case of nodal hypersurfaces. A similar result for hypersurfaces $D$ having arbitrary isolated singularities is obtained in \cite{D3}.
 
Illustrations of how to apply Theorem \ref{onenode} are given in Example \ref{exonenode}.
Using Theorems \ref{onenode} and \ref{thm5} and Corollary \ref{corA1}, we get the following information on the {\it position of the nodes of a nodal curve}.

\begin{cor}
\label{corC1}
Let $C:f=0$ be a degree $N$ nodal curve in $\PP^2$ and let $\NN$ denote the set of its nodes. Then
$\defect S_k(\NN)=0$ for $k>N-3$ and $\defect S_{N-3}(\NN)=r-1$, where $r$ is the number of irreducible components of $C$.

Moreover, if the curve $C$ is in addition rational, all the defects $\defect S_k(\NN)$ are completely determined by the degree $N$ and the number of nodes $n(C)$.
\end{cor}

In fact, a recent result by Kloosterman, see Proposition 3.6 in \cite{RK}, implies that the first part of Corollary \ref{corC1} holds for any curve $C$ with the property that any singular point of $C$ which is not a node is a unibranch singularity, see Remark \ref{RK} for more details on this.

\medskip

In the last section we use Theorem \ref{onenode} to determine the pole order filtration $P^*$ on the cohomology groups $H^*(U)$ and the corresponding spectral sequences when $D$ is a nodal surface. In particular, we get the following.

\begin{thm}
\label{thm6} Let $S:f=0$ be a nodal surface in $\PP^3$ of degree $N$ and let $\NN$ denote the set of its nodes. Then, if $U=\PP^3 \setminus S$ and $S_s$ is a smooth surface of degree $N$ in $\PP^3$, the following hold.
$$\dim Gr_P^2(H^3(U))=h^{1,1}(S_s)-1-\defect S_{N-4}(\NN)$$
and 
$$\dim Gr_F^2(H^3(U))=h^{1,1}(S_s)-1-|\NN|.$$
In particular, $P^2H^3(U)=F^2H^3(U)$ if and only if the nodal surface $S$ is smooth or $N<4$.
\end{thm}
This result complements the results in \cite{DSW} (where arbitrary dimensions are considered but only degrees $N=3$ and $N=4$), in the case of nodal surfaces, and answers the question asked there whether the inequality $P^2H^3(U) \ne F^2H^3(U)$ holds for any surface with $|\NN|=1$ and $N \geq 4$.

\bigskip

Numerical experiments with the CoCoA package \cite{Co} and the Singular package \cite{Sing} have played a key role in the completion of this work.

\section{Pole order filtrations, spectral sequences and Koszul complexes} \label{sec:two}

Let $X$ be a smooth complex quasi-projective variety and  $D \subset X$ a reduced divisor. We denote by $i: D \to X$ and $j:U \to X$ the corresponding inclusions, where $U=X \setminus D$. Let $\Omega_X^*$ (resp. $\Omega_U^*$)  denote the de Rham sheaf complex of regular differential forms on $X$ (resp. on $U$). Then Grothendieck's Theorem says that
\begin{equation} 
\label{GT}
\HH^*(U,\Omega_U^*)=H^*(U)
\end{equation}
where $\C$-coefficients are used for the cohomology groups in this paper unless indicated otherwise.
Moreover, as explained in \cite{DD}, the isomorphism $j_*\Omega_U^*=Rj_*\Omega_U^*$, which is due to the fact that $j$ is an affine morphism, implies a natural identification
\begin{equation} 
\label{GT1}
\HH^*(X,j_*\Omega_U^*)=H^*(U).
\end{equation}
The sheaf complex $j_*\Omega_U^*$ has a natural decreasing filtration, called the {\it pole order filtration}, given by
$P^sj_*\Omega_U^p=0$ if $p<s$ and
\begin{equation} 
\label{po}
P^sj_*\Omega_U^p=\Omega_X^p((p-s+1)D)
\end{equation}
if $p \geq s$, see \cite{DD}. In other words, a rational differential form $\omega$ is in $P^sj_*\Omega_U^p$
if it has a pole of order at most $(p-s+1)$ along the divisor $D$ (with a special attention needed for the case of $p=s-1$.) A word of warning: the corresponding filtration is denoted by $F$ in \cite{Dc}
and is slightly different. However the proof of the main results from \cite{Dc} or \cite{D1} quoted below apply word for word to the present setup.

Using the filtration \eqref{po}, we define the pole order filtration on the cohomology of $U$
by setting
\begin{equation} 
\label{pf}
P^sH^*(U)=\im(\HH^*(X,P^sj_*\Omega_U^*) \to \HH^*(X,j_*\Omega_U^*)=H^*(U)).
\end{equation}
The main result from \cite{DD} is the following. See also \cite{Sa1} for another proof and conditions
for equality.

\begin{thm}
\label{thmDD}
Assume that the smooth variety $X$ is proper and let $F$ denote the Hodge filtration on the cohomology of $U$. Then $F^sH^*(U) \subset P^sH^*(U)$ for any $s$.
\end{thm}
From now on consider the case $X=\PP^n$ and recall that Bott's vanishing theorem gives us
\begin{equation} 
\label{Bott}
H^k(X,\Omega^p_X(sD))=0
\end{equation}
for any $k>0$, $s>0$, see \cite{Bo}.
The polar filtration, even if it is an infinite filtration, gives rise to a spectral sequence
\begin{equation} 
\label{sp1}
E^{p,q}_1(U)=\HH^{p+q}(X,Gr_P^p(j_*\Omega_U^*))
\end{equation}
whose limit term is exactly
\begin{equation} 
\label{sp2}
E^{p,q}_{\infty}(U)=Gr_P^p(H^{p+q}(U)).
\end{equation}
Using now the standard spectral sequence
\begin{equation} 
\label{sp3}
E^{p,q}_1=H^p(X,Gr_P^s(j_*\Omega_U^q)) \Rightarrow \HH^{p+q}(X,Gr_P^s(j_*\Omega_U^*))
\end{equation}
and the vanishings implied by \eqref{Bott}, we get a description of the $E_1$-term of our spectral sequence without involving
hypercohomology groups, namely
\begin{equation} 
\label{sp4}
E^{p,q}_1(U)=H^{p+q}(H^0(X,Gr_P^p(j_*\Omega_U^*))).
\end{equation}
This expression for $E^{p,q}_1(U)$ can be interpreted as follows. Let $A^*(U)=H^0(X,j_*\Omega_U^*)$
be the de Rham complex of regular forms defined on the affine open set $U$. It follows from Grothendieck's Theorem \ref{GT}, that one has
\begin{equation} 
\label{GT2}
H^m(A^*(U))=H^m(U)
\end{equation}
for any integer $m$. On the other hand, we have a very explicit description of these rational differential forms defined on $U$. Let $f=0$ be a reduced equation for the divisor $D$ and let $N$ be the degree of the homogeneous polynomial $f$. Denote by $\Omega^p=H^0(\C^{n+1}, \Omega_{\C^{n+1}}^p)$ the global (polynomial) differential $p$-forms on $\C^{n+1}$, regarded as a graded $S$-module in the usual way (i.e. $\deg(hdx_{i_1}\wedge...\wedge dx_{i_q}) =p+q$ if $h \in S_p$).
Then a differential $p$-form  $ \omega \in A^p(U)$, for $p \geq 0$, is given by 
\begin{equation} 
\label{diff1}
\omega=\frac{\Delta(\gamma)}{f^s}
\end{equation}
for some integer $s>0$, $\gamma \in \Omega^{p+1}_{sN}$ and $\Delta: \Omega^{p+1} \to \Omega^{p}$ the $S$-linear map given by the contraction with the Euler field, see Chapter 6 in \cite{D1} for details.
When  $\omega$ is not a constant function on $U$, case covered by $s=1$ and $\gamma= a \cdot df$ for $a \in \C$, the minimal $s$ in this formula is by definition the order of $\omega$ along the divisor $D$.
We can define a polar filtration on the complex $A^*(U)$ by setting
$P^sA^p(U)=0$ if $p<s$ and
\begin{equation} 
\label{po2}
P^sA^p(U)=\{ \omega=\frac{\Delta(\gamma)}{f^{p-s+1}}~~|~~ \gamma \in \Omega^{p+1}_{(p-s+1)N} \}    
\end{equation}
if $p \geq s$. This decreasing filtration induces a spectral sequence
\begin{equation} 
\label{sp5}
E^{p,q}_1(A)=H^{p+q}(Gr_P^p(A^*(U))).
\end{equation}
Using Bott's vanishing \ref{Bott} and  \eqref{sp4}, we see that this new spectral sequence 
coincides with the spectral sequence $E^{p,q}_1(U)$, in particular they induce both the same filtration on their common limit which is $H^*(U)$.

Note that $A^0(U)$ (resp. $E^{0,0}_1(A)=H^{0}(Gr_P^0(A^*(U)))$) contains the constant functions on $U$. Let us denote by $\tilde A^*(U)$
(resp. $E^{p,q}_1(\tilde A)$) the complex (resp. the spectral sequence) obtained from the above complex
$A^*(U)$ (resp. spectral sequence $E^{p,q}_1(A)$)
by replacing $A^0(U)$ (resp. $E^{0,0}_1(A)$) by $A^0(U)/\C$ (resp. $E^{0,0}_1(\tilde A)=E^{0,0}_1(A)/\C$). It is clear that the cohomology of the complex $\tilde A^*(U)$ (resp. the limit of the spectral sequence $E^{p,q}_r(\tilde A)$) is $\tilde H^*(U)$, the reduced cohomology of $U$.

\bigskip

It turns out that the $E_1$-term of the spectral sequence $E^{p,q}_r(\tilde A)$ can be described in terms of the {\it Koszul complex} of the partial derivatives $f_j$ of $f$ with respect to the variable $x_j$ for $j=0,...,n$,
see \cite{Dc}, \cite{D1}, Chapter 6 and Remark 2.10 in \cite{DS1}.

This Kozsul complex can be represented by the complex of graded $S$-modules
\begin{equation} 
\label{kc}
K^*(f): 0 \to \Omega^0 \to \Omega^1 \to ... \to \Omega^{n+1} \to 0
\end{equation}
where the differentials  are given by the wedge product with the differential $df$ and hence these differentials are homogeneous of degree $N$.
This complex has a natural subcomplex
\begin{equation} 
\label{kc1}
K'^*(f): 0 \to \Omega'^0 \to \Omega'^1 \to ... \to \Omega'^{n+1} \to 0
\end{equation}
where $\Omega'^p=\oplus_{k \geq 0}\Omega^p_{kN}.$

\medskip

Consider the associated double complex $(B,d',d'')$, with $B^{s,t}=\Omega^{s+t+1}_{(t+1)N}$ for $t\geq 0$ and $-1 \leq s+t \leq n$ and $B^{s,t}=0$ otherwise, and differentials $d'=d$, the exterior derivative of a form, and $d''(\omega)=-df\wedge \omega.$
Note that $d'd''+d''d'=0$ and let $(B^*, D_f=d'+d'')$ be the associated total complex of this double complex. In fact the complex $B^*$ is the same as the reduced version of the subcomplex $K'^*$, but with a new differential.

As for any total complex, it comes with two natural decreasing filtrations, one of them being
$$F^pB^k=\oplus_{s \geq p-1}B^{s,k-s}.$$
The contraction operator $\Delta$ defines a morphism of filtered complexes $\delta: B^* \to \tilde A^*(U)$ by setting
\begin{equation} 
\label{delta}
\delta(\omega)=\frac{\Delta(\omega)}{f^{t+1}} \text{   for   } \omega \in B^{s,t}.
\end{equation}
With this notation, we have the following result, \cite{Dc} \cite{D1}, Chapter 6 and Remark 2.10 in \cite{DS1}.

\begin{prop}
\label{prop1}
Let $E_r^{p,q}(f)$ be the $E_1$-spectral sequence associated to the filtration $F$ on $(B^*,D_f)$.
Then the following hold.

\noindent (i) The morphism $\delta$ induces an isomorphism of $E_1$-spectral sequences
$$E_r^{p,q}(f) \to  E_r^{p,q}(\tilde A).$$
\noindent (ii) There is a natural identification 
$$E_1^{p,q}(f)=H^{p+q+1}(K^*(f))_{(q+1)N}.$$

\end{prop}

\begin{rk}
\label{rkA}

\noindent (i)
In the case $X=\PP^n$, it is known that $F^1H^k(U)=H^k(U)$ for any integer $k>0$, see Theorem (2.2)
in \cite{Dc} (there is a sign $=$ missing in the statement, but the proof of the equality is clearly done) or the proof of Corollary 1.32 on pp. 185-186 in \cite{D1}.

\noindent (ii) One has $P^{k+1}H^k(U)=0$ for any integer $k>0$. To see this, just use the fact that the hypercohomology of a sheaf complex $\F^*$ with $\F^j=0$ for $j<p$ satisfies $\HH^j(\F*)=0$ for $j<p$. In particular $P^2H^1(U)=0$, i.e. we always have $Gr_P^1(H^{1}(U))=H^1(U)$ and $Gr_P^j(H^{1}(U))=0$ for $j \ne 1$.
\end{rk}

Assume now that the hypersurface $D$ has only isolated singularities. The nonzero terms in the $E_1$-term of the spectral sequence $E_r^{p,q}(f)$ are sitting on two lines, given by $L: p+q=n$ and $L': p+q=n-1$.
Indeed, one has to use the fact that in this case $H^m(K^*(f))=0$ for $m<n$, see \cite{Gr1}, \cite{Sai}.

For a term $E_1^{p,q}(f)$ situated on the line $L$, we have
$$E_1^{p,q}(f)=H^{n+1}(K^*(f))_{(q+1)N}=M(f)_{(q+1)N-n-1}.$$

We describe now the terms on the line $L'$. In order to do this, let $f_s \in S_N$ denote a polynomial of degree $N$ defining a smooth hypersurface in $\PP^n$. 

It is easy to show that
\begin{equation} 
\label{eq22}
t^NHP(H^n(K^*(f)))(t)= HP(H^{n+1}(K^*(f)))(t)-HP(H^{n+1}(K^*(f_s)))(t),
\end{equation}
using the fact that Euler characteristics do not change when replacing a (finite type) complex by its cohomology.
Note also that
\begin{equation} 
\label{eq23}
HP(H^{n+1}(K^*(f_s)))=t^{n+1}HP(M(f_s))=t^{n+1} \cdot \frac{(1-t^{N-1})^{n+1}}{(1-t)^{n+1}}
\end{equation}
is completely determined by the degree $N$.

\medskip

It follows that the term $E_1^{p,q}(f)=H^{n}(K^*(f))_{(q+1)N}$ situated on the line $L'$ has dimension
\begin{equation} 
\label{eq24}
\dim H^{n}(K^*(f))_{(q+1)N}= \dim  M(f)_{(q+2)N-n-1} - \dim  M(f_s)_{(q+2)N-n-1}.
\end{equation}

\medskip

We want now to relate the spectral sequence $E_r^{p,q}(A)$ to some simpler, locally computable spectral sequences in the case when $D$ has only isolated singularities, say at the points $a_1,...,a_m$. Consider the morphism of restriction
$$\rho: Gr_P^p(j_*\Omega_U^*) \to i_{1*}Gr_P^p((j_*\Omega_U^*)/\Omega_X^*)| \Sigma$$
obtained by factoring out the regular forms, then taking the restriction from $X$ to the singular locus $\Sigma$ of $D$, and then extending via $i_{1*}$, where $i_{1}:\Sigma \to X$ is the inclusion.
For $p<0$ this morphism is easily seen to be a quasi-isomorphism, i.e. it induces isomorphisms at stalk level. For $p=0$, the kernel $K_{\rho}$ of $\rho$ is the sheaf $\Omega^0_X=\OO_{X}$ (placed in degree zero). 
We know that, in the case $X=\PP^n$,
$$\HH^q(X, \OO_{X})=H^q(X,\OO_{X})=0$$
for $q>0$. 
It follows that the morphisms 
$$\rho^k: \HH^k(X, Gr_P^0(j_*\Omega_U^*)) \to \HH^k(X, i_{1*}Gr_P^0((j_*\Omega_U^*/\Omega_X^*)|\Sigma)$$
are isomorphisms for any $k \geq 1$.

As explained in \cite{Dc} (with the notable difference that in loc.cit. there is no quotient taken, which leads to an  infinite dimensional $E_1$-term), the complex $((j_*\Omega_U^*)/\Omega_X^*)| \Sigma$ is the direct sum of the complexes 
$\tilde A^*(D,a_j)$ for $j=1,...,m$, where each $\tilde A^*(D,a_j)$ is the local analog of the complex $\tilde A^*(U)$ above.
These complexes come with a pole order filtration defined exactly as in the global case, and for each $j$ there is an $E_1$-spectral sequence $E_r(D,a_j)$ with 
$$E_1^{p,q}(D,a_j)=H^{p+q}(Gr_P^p(\tilde A^*(D,a_j)))$$
and converging to $\tilde H^*(B_j \setminus D)$, where $B_j$ is a small ball in $X$ centered at $a_j$.

It follows that $\rho$ induces a morphism of $E_1$-spectral sequences
$$\rho^{p,q}:E_1^{p,q}(A) \to \oplus_{j=1,m}E_1^{p,q}(D,a_j)$$
with the property that $\rho^{p,q}$ is an isomorphism for any $p\leq 0$ and $p+q \geq 1$.

Moreover, when each singularity $(D,a_j)$ is weighted homogeneous, it follows from the description of the local spectral sequence $E_1^{p,q}(D,a_j)$, see Example 3.6 in \cite{Dc}, that all the differentials $d_1:E_1^{n-1-t,t}(D,a_j) \to E_1^{n-t,t}(D,a_j)$
are isomorphism for $t \geq n-1$.

In this way we have proved the following improvement of Theorem (3.9) in \cite{Dc}. (For the converse
claim in $(iii)$ see Corollary (3.10) in \cite{Dc}).

\begin{thm}
\label{thm2}

\noindent (i) Let $D$ be a hypersurface in $\PP^n$ for $n\geq 2$, having only isolated singularities. Then morphism of $E_1$-spectral sequences
$$\rho^{p,q}:E_1^{p,q}(A) \to \oplus_{j=1,s}E_1^{p,q}(D,a_j)$$
is an isomorphism for any $p\leq 0$ and $p+q \geq 1$.

\noindent (ii) If in addition the  singularities of $D$ are weighted homogeneous,
then in the spectral sequence $E_1^{p,q}(A)$ the differential 
$$d_1:E_1^{n-1-t,t}(A) \to E_1^{n-t,t}(A)$$
is injective for $t=n-1$ and is bijective for $t \geq n.$

\noindent (iii) If $D$ is a reduced curve in $\PP^2$, then $D$ has only isolated weighted homogeneous singularities if and only if the $E_1$-spectral sequences $E^{p,q}_r(U)$, $E^{p,q}_r(\tilde A)$  and $E^{p,q}_r(f)$
degenerate at the $E_2$-term, i.e. $E_2=E_{\infty}$ for any of these $E_1$-spectral sequences.

\end{thm}

This result, especially the parts (ii) and (iii), is perhaps related to the results in \cite{CN}
and \cite{CMNJ}.

\section{Some examples of spectral sequences in the case of plane curves} \label{sec3}

Let $C:f=0$ be a reduced curve in $\PP^2$ of degree $N$.
Let $C_j:f_j=0$ for $j=1,...,r$ be the irreducible components of $C$. The complement $U$ has at most three non-zero cohomology groups.
The first of them, $H^0(U)$ is $1$-dimensional and of Hodge type $(0,0)$, so nothing interesting here.
Moreover $\tilde H^0(U)=0$.

The second one, $H^1(U)$ is $(r-1)$-dimensional and, for $r>1$, is of Hodge type $(1,1)$ by Remark \ref{rkA}. It follows that in this case $P^1H^1(U)=F^1H^1(U)=H^1(U)$.
Moreover, $H^1(U)$ has a basis given by 
\begin{equation} 
\label{basis1}
\omega _j=\frac{df_j}{N_jf_j}-\frac{df_r}{N_rf_r}
\end{equation} 
for $j=1,...,r-1$, where $N_j=\deg(f_j)$, see \cite{Dc}, Example (4.1).

\begin{ex}
\label{FP0}

We discuss first the case when $C:f=0$ is a nodal curve in $\PP^2$ of degree $N$.
Using Corollary (0.12) in \cite{Sa1} for $X=\PP^2$, $i=2$ it follows that $P^2H^2(U)=F^2H^2(U)$,
since for a nodal curve $\al_f=1$.
Now we look at the nonzero terms in the $E_1$-term of the spectral sequence $E_r^{p,q}(f)$.
They are sitting on two lines, given by $L: p+q=2$ and $L': p+q=1$.

\medskip

We look first at the terms on the line $L$. The term $E_1^{2,0}(f)=H^{3}(K^*)_{N}$ is isomorphic as a $\C$-vector space to $M(f)_{N-3}$, hence has dimension $g$, as defined in \eqref{dim1},
which is determined by $N=\deg (f)$ alone.
Hence, the corresponding limit term $E^{2,0}_{\infty}(U)=P^2H^{2}(U)$ has the  dimension at most $g$. On the other hand, $\dim F^2H^{2}(U)=g$ , see Theorem 2.2 in \cite{DSW} or a direct proof in Proposition 4.1 in \cite{DSt}. This gives an alternative proof of the equality $F^2H^{2}(U)=P^2H^{2}(U)$ in this case.

\medskip

The term $E_1^{1,1}(f)=H^{3}(K^*)_{2N}$ is isomorphic to $M(f)_{2N-3}$. To compute its dimension, note that we have $\dim (I/J_f)_{2N-3}=Gr_F^1(H^{2}(U))$ by Theorem 2.2 in \cite{DSW}, where $I$ is the ideal in $S$ of polynomials vanishing at all the singular points of $C$.
It was shown in Proposition 4.1 in \cite{DSt} that
$$\dim (I/J_f)_{2N-3}=\sum_{j=1,r}g_j$$
where $g_j$ is the genus of the normalization of the curve $C_j$, for $j=1,...,r$.
On the other hand, we showed in Lemma 4.2 in \cite{DSt} that 
$\dim(S/I)_{2N-3}=n(C),$
the total number of nodes of $C$.
It follows that
\begin{equation} 
\label{dim2}
\dim M(f)_{2N-3}=n(C)+\sum_{j=1,r}g_j.
\end{equation}
Moreover, the corresponding limit term $E^{1,1}_{\infty}(U)=Gr_F^1(H^{2}(U))= Gr_P^1(H^{2}(U)) $ has  dimension $\sum_{j=1,r}g_j$
as noted above.

\medskip

The  term $E_1^{2-q,q}(f)=H^{3}(K^*)_{(q+1)N}$ for $q\geq 2$ is isomorphic to $M(f)_{(q+1)N-3}$, which has dimension $n(C)$.
Moreover, the corresponding limit terms $E^{2-q,q}_{\infty}(U)=Gr_F^{2-q}(H^{2}(U))$ clearly vanish for $q\geq 2$.

\medskip

We look now at the terms on the line $L'$. It follows that the term $E_1^{1,0}(f)=H^{2}(K^*)_{N}$ has dimension
$n(C)+\sum_{j=1,r}g_j-g$,
via \eqref{dim2},  \eqref{dim1} and the duality $\dim  M(f_s)_{2N-3}=\dim  M(f_s)_{N-3}$.
If we compare with the proof of Proposition 4.1 in \cite{DSt}, we see that the total number of nodes $n(C)$ is given by $\sum_{j=1,r}n_j+ \sum_{1\leq i<j \leq r}d_id_j$ where $n_j$ is the number of nodes on the curve $C_j$ and $d_k$ is the degree of the curve $C_k$. Using the formula (4.1) in the proof of Proposition 4.1 in \cite{DSt} and Remark \ref{rkA}, we conclude that
\begin{equation} 
\label{dim3}
\dim E_1^{1,0}(f)=\dim E_{\infty}^{1,0}(f)=r-1.
\end{equation}

\medskip

The dimension of the other terms $E_1^{1-q,q}(f)=H^{2}(K^*)_{(q+1)N}$ for $q\geq 1$ is equal to  $n(C)$.
Moreover, the corresponding limit terms $E^{1-q,q}_{\infty}(U)=Gr_F^{1-q}(H^{1}(U))$ clearly vanish for $q\geq 1$.

It follows that the differential $d_1: E_1^{1,0}(f) \to E_1^{2,0}(f)$ is the zero map
(not to decrease the dimension of $E_2^{1,0}(f)$, which is the dimension of the limit),
a fact which is not shared by curves with general weighted homogeneous singularities as seen in Examples
\ref{FP1} and \ref{FP2} below.
The other differentials $d_1: E_1^{1-q,q}(f) \to E_1^{2-q,q}(f)$ for $q\geq 1$ are all injective
(any nonzero kernel would kill some terms needed in the limit via some $d_r$ with $r \geq 2$), as it happens for any curve with weighted homogeneous singularities in view of Theorem \ref{thm2}.

\end{ex}

\begin{ex}
\label{FP1} Consider the curve $C:x(x^2y+xy^2+z^3)=0$, which is the union of a smooth cubic $C: x^2y+xy^2+z^3=0 $ and an inflectional tangent $L:x=0$.
Then it is easy to see that $\tilde H^0(U)=0$, $H^1(U)$ is $1$-dimensional, and $H^2(U)$ is $2$-dimensional,
with classes of Hodge type $(2,1)$ and $(1,2)$. In particular, $F^2H^2(U)$ is $1$-dimensional.

On the other hand, the spectral sequence $E_1(f)$ has the following nonzero terms:
$E_1^{1,0}$ which is $2$-dimensional, $E_1^{2,0}$ which is $3$-dimensional, and all $E_1^{p,q}$ for $p+q=1$ or $p+q=2$ and $q>0$, which are $5$-dimensional, since $\tau(C)=5$. The computation for the other dimensions are based on formula \eqref{eq23} and a computation using CoCoA \cite{Co} or Singular \cite{Sing} of the Hilbert-Poincar\'e series
$$HP(M(f))(t)=1+3t+6t^2+7t^3+6t^4+5t^5+...$$
with stabilization threshold $st(C)=5$.
It follows that $d_1: E_1^{1,0} \to E_1^{2,0}$ has a $1$-dimensional kernel
$E_2^{1,0}=H^1(U)$, and a $2$-dimensional cokernel $E_2^{2,0}=P^2H^2(U)$. In particular, the inclusion
$F^2 \subset P^2$ is strict on $H^2(U)$ as mentioned in \cite{Dc}, Remark (2.6).
\end{ex}

\begin{ex}
\label{FP2}

 Consider now the irreducible curve $C:x^2y^2+xz^3+yz^3=0$, which has two cusps $A_2$ as singularities.
Then it is easy to see that $\tilde H^0(U)=0=H^1(U)$, and $H^2(U)$ is $2$-dimensional,
with classes of Hodge type $(2,1)$ and $(1,2)$. In particular, $F^2H^2(U)$ is $1$-dimensional.

On the other hand, the spectral sequence $E_1(f)$ has the following nonzero terms:
$E_1^{1,0}$ which is $1$-dimensional, $E_1^{2,0}$ which is $3$-dimensional, and all $E_1^{p,q}$ for $p+q=1$ or $p+q=2$ and $q>0$, which are $4$-dimensional, since $\tau(C)=4$. Indeed, the computation using CoCoA \cite{Co} or Singular \cite{Sing} yields in this case
$$HP(M(f))(t)=1+3t+6t^2+7t^3+6t^4+4t^5+...$$
with stabilization threshold $st(C)=5$.
It follows that $d_1: E_1^{1,0} \to E_1^{2,0}$ is injective
and has a $2$-dimensional cokernel $E_2^{2,0}=P^2H^2(U)$. In particular, the inclusion
$F^2 \subset P^2$ is strict on $H^2(U)$ as mentioned in \cite{DS1}, Remark (2.5).
\end{ex}

\begin{ex}
\label{FP3}

 Consider now the irreducible curve $C:x^3z^4+xy^5z+x^7+y^7=0$, which has a {\it non weighted homogeneous singularity}
located at $(0:0:1)$ with Milnor number $\mu= 12$ and Tjurina number $\tau=11$.
Then it is easy to see that $\tilde H^0(U)=0=H^1(U)$, and $H^2(U)$ has dimension $18$, with classes of Hodge type $(2,1)$ and $(1,2)$. In particular $\dim F^2H^2(U)=9$.

On the other hand, the spectral sequence $E_1(f)$ has the following nonzero terms:
 $E_1^{2,0}$ which is $15$-dimensional, and all $E_1^{p,q}$ for $p+q=1$ or $p+q=2$ and $q>0$, which are $11$-dimensional, since $\tau(C)=11$, except $E_1^{1,1}$ which is again $15$-dimensional. Indeed, the computation using CoCoA \cite{Co} or Singular \cite{Sing} yields in this case
$$HP(M(f))(t)=1+3t+6t^2+10t^3+15t^4+21t^5+25t^6+27t^7+$$
$$+27t^8+25t^9+21t^{10}+15 t^{11}+12t^{12}+11t^{13}...$$
with stabilization threshold $st(C)=13$.
It follows that $d_1: 0=E_1^{1,0} \to E_1^{2,0}$ is the zero map, hence
$\dim E_2^{2,0}=15$. The other differentials $d_1: E_1^{1-t,t} \to E_1^{2-t,t}$
for $t \geq 1$ have a  $1$-dimensional kernel, use Theorem \ref{thm2}, part (i), and Proposition (3.4), Example (3.5) (i) and Corollary (4.3) in  \cite{DProc}, where it is shown that in this case the differentials $d_2: E_2^{1-t,t} \to E_2^{3-t,t-1}$ are injective for $t>0$ in the local setting. It follows that $E_3=E_{\infty}$ has the following nonzero terms: $E_3^{2,0}$ of dimension $14$, and  $E_3^{1,1}$ of dimension $4$. In particular, one has 
$$\dim F^2H^2(U)=9<14=\dim P^2H^2(U).$$

\end{ex}

\section{The syzygies of  nodal hypersurfaces} \label{sec4}

First we give a geometric interpretation of a syzygy  $R_m$ as in \eqref{rel1} in the case $n=2$ using section (2.1) in
\cite{BDS}. Let $F_f$ be the Milnor fiber of $f$, which is the smooth affine surface in $\C^3$ given by the equation $f(x,y,z)=1$. Then there is a monodromy isomorphism $h:F_f \to F_f$ given by multiplication by $\lambda=\exp(2\pi i/N)$ and an induced monodromy operator $h^1:H^1(F_f) \to H^1(F_f)$.
The eigenvalues of $h^1$ are exactly the $N$-th roots of unity and for each $k=0,1,...,N-1$ there is a rank one local system $L_k$ on $U$ such that
\begin{equation} 
\label{rel2}
H^*(F_f)_{\lambda^k}=H^*(U,L_k)
\end{equation}
where in the LHS we have the corresponding eigenspace and in the RHS we have the twisted cohomology of $U$ with coefficients in $L_k$, see for details \cite{D2}, Proposition 6.4.6.

Let $\LL_k$ be the Deligne extension of $L_k$ over the nodal curve $C$ such that the eigenvalues of the residue of the connection are contained in the interval $[0,1)$. In our case, the line bundle $\LL_k$ is precisely
$\OO_{\PP^2}(-k)$, see (2.1.2) in \cite{BDS} and we have the following relation with the Hodge filtration on $H^*(F_f)$:
\begin{equation} 
\label{rel3}
Gr^p_FH^{p+q}(F_f)_{\lambda^k}=H^q(\PP^2,\Omega^p_{\PP^2}(logC)\otimes \LL_k)
\end{equation}
see (2.1.1) in \cite{BDS}. In particular, we get
\begin{equation} 
\label{rel4}
Gr^1_FH^{1}(F_f)_{\lambda^k}=H^0(\PP^2,\Omega^1_{\PP^2}(logC)\otimes \LL_k).
\end{equation}
Now the curve $C$ being nodal, it follows that $H^{1}(F_f)_{\lambda^k}=0$ for $k=1,...,N-1$, see Corollary 6.4.14 in \cite{D2} for a stronger result.

Assume now that we have a nonzero syzygy  $R_m$ as in \eqref{rel1} with $m<N-2$. Consider the nonzero $2$-form $\omega \in \Omega^2_{m+2}$  given by $\omega=ady\wedge dz-bdx \wedge dz+c dx \wedge dy$ and note that $df \wedge \omega=0$. The $1$-form 
\begin{equation} 
\label{rel5}
\al=\frac{\Delta(\omega)}{f}
\end{equation}
is an element of $H^0(\PP^2,\Omega^1_{\PP^2}(logC)\otimes \LL_k)$, with $k=N-2-m>0$.
To see this, use the formula for $d\al$ given in (1.10), p. 181 in \cite{D1}. Moreover, $\al \ne 0$ since  the kernel of $\Delta: \Omega^2 \to \Omega^1$ is the free $S$-module spanned by 
$\sigma=\Delta(dx \wedge dy \wedge dz)$ and $df \wedge \sigma=Nfdx \wedge dy \wedge dz\ne 0$.
But this is in contradiction to $H^{1}(F_f)_{\lambda^k}=0$ in view of \eqref{rel4}.

\bigskip

Next we'll describe all the syzygies  $R_m$ as in \eqref{rel1} with $n=2$ and $m=N-2$. This is the same as
describing $H^2(K^*(f))_N$, and we know from the previous section that $\dim H^2(K^*(f))_N=r-1$, see \eqref{dim3}. This means essentially to lift the basis $\omega_j$ in \eqref{basis1} to a basis of
$H^2(K^*(f))_N$. Note that
$$\omega_j=\frac{\al_j}{N_jN_rf}$$
where $\al_j=N_rf_1...\hat f_j...f_rdf_j-N_jf_1....\hat f_rdf_r$ for $j=1,...,r-1$ and $\hat f_j$
means that the factor $f_j$ is missing.
Define $\be_j=-f_1...\hat f_j....\hat f_r df_j \wedge df_r$ and note that
$$\Delta(\be_j)=-f_1...\hat f_j....\hat f_r \Delta(df_j \wedge df_r)=\al_j.$$
For $r=2$, $\be_1$ is a good lifting since $df \wedge \be_1=0$ and we are done.
However, for $r>2$, $\be_j$ is not a good lifting, since in general one has 
$$df \wedge \be_j=-\sum_{k \ne j;~~k\ne r}f^2/(f_kf_jf_r)df_k\wedge df_j \wedge df_r=fg_jdx \wedge dy \wedge dz$$
for some $g_j \in S_{N-3}$ which nonzero in general. (A formula for $g_j$ is given in Theorem \ref{thm5} below 
using the Jacobian determinant $Jac(f_k,f_j,f_r)$ of the three functions $f_k,f_j,f_r$ with respect to $x,y$ and $z$).

To correct this problem, we look for a modification of the form
$$\gamma_j=\be_j+h_j\sigma$$
where $h_j \in S_{N-3}$ and $\sigma=\Delta(dx \wedge dy \wedge dz)$ as above.
Now $df \wedge \gamma_j=(fg_j+Nfh_j)dx \wedge dy \wedge dz=0$ if we choose $h_j=-g_j/N$.
The resulting $\gamma_j$ for $j=1,...,r-1$ yield a basis of $H^2(K^*(f))_N$.

Hence we have proved the following result. 

\begin{thm}
\label{thm5}
Let $C:f=0$  be a nodal curve of degree $N$ in $\PP^2$.
Then $H^2(K^*(f))_q=0$ for any $q<N$, $H^2(K^*(f))_N$ is $(r-1)$-dimensional and a basis for it is given by
$$\gamma_j=-f_1...\hat f_j....\hat f_r df_j \wedge df_r+h_j \sigma,$$
for $j=1,...,r-1$.
Here $r$ is the number of irreducible components of $C$, $f_j=0$ are reduced equations for these components,
$\sigma=\Delta(dx \wedge dy \wedge dz)$, $h_1=0$ if $r=2$ and
$$h_j=\frac{\sum_{k\ne j; ~~k \ne r}f/(f_kf_jf_r)Jac(f_k,f_j,f_r)}{N}$$
if $r>2$.

\end{thm}
For an arbitrary curve $C$ having $r$ irreducible components $C_j:f_j=0$, the above elements $\gamma_j$
yield $r-1$ linearly independent elements in $H^2(K^*(f))_N$, which are killed by $d_1$.
It may happen that $\dim H^2(K^*(f))_N>r-1$, as we have seen in Example \ref{FP1}.

\medskip
The corresponding vanishing result in the general case of nodal hypersurfaces is considered in \cite{DSt3}, but in this general case there is no description of an explicit basis of the lowest degree (possibly nonzero) syzygies as in Theorem \ref{thm5} above.
For an alternative proof of the vanishing part (without using Hodge theory)
in a more general curve setting, see Eisenbud and Ulrich \cite{EU}.

\begin{ex}
\label{lowdegrel}
In this example we look at some curves having low degree relations $R_m$ as in \eqref{rel1}.

\noindent (i) It is clear that a curve $C:f=0$ admits a relation of degree $m=0$ if and only if up-to a linear coordinate change we have that the equation $f$ is independent of $z$. In this case
$$HP(M(f))(t)=\frac{(1-t^{N-1})^2}{(1-t)^3}.$$
Hence $ct(C)=N-2$ (this is the minimal possible value) and $st(C)=2N-4.$

\noindent (ii) The curve $C: x^py^q+z^N=0$ for $p+q=N$ admits an obvious relation of degree one, namely
$$qxf_x-pyf_y=0.$$
In this case $ct(C)=N-1$.

\noindent (iii) The curve $C: z^p(x^q+y^q)+x^N+y^N=0$ for $p+q=N$ admits an obvious relation of degree $2p$,
namely
$$z^{p-1}x(qz^p+Ny^p)f_x+ z^{p-1}y(qz^p+Nx^p)f_y-\frac{1}{p} (qz^p+Ny^p)(qz^p+Nx^p)f_z=0  .$$
It is easy to see that this relation is not a consequence of the trivial relations $T_{ij}$ in \eqref{tij}. On the other hand, a computation in the case $N=7$, $p=4$ shows that
$$HP(M(f))(t)=1+3t+6t^2+10t^3+15t^4+21t^5+25t^6+27t^7+$$
$$+27t^8+25t^9+21t^{10}+16 t^{11}+12t^{12}+9t^{13}+8t^{14}+...$$
with stabilization threshold $st(C)=14$. It follows that $ct(C)=10$, which implies via \eqref{REL} that
$mdr(C)=5$, i.e. the above relation has not minimal degree in general.
However this is the case for $p=N-2$, when the curve $C$ has a node at $(0:0:1)$ and the corresponding relation has degree $2N-4$.

\end{ex}

This is a very special case of Theorem \ref{onenode} stated in the Introduction, and which we prove now.

\proof

Choose the coordinates on $\PP^n$ such that $H_0:x_0=0$ is transverse to $D$, i.e. the intersection $H_0 \cap D$ is smooth. It follows as in \cite{CD} that the partial derivatives $f_1$,...,$f_n$ of $f$ form a regular system in $S$, in particular they vanish at a finite set of points on $\PP^n$, say $p_1,...,p_r$. A part of these points, say $p_j$ for $j=1,...,q$ are  the nodes on $D$, i.e. the points in the set $\NN$. It follows that the divisors $D_j:f_j=0$ for $j=1,...,n$ intersect transversely at any point $p_j \in \NN$. To see this, one may work in the affine chart $x_0=1$, where $x_1,...,x_n$ may be used as coordinates and use the definition of nodes as the singularities where the hessian of a (local) equation is nonzero.

Assume we have a nonzero element in $H^{n}(K^*(f))_{nN-n-1-k}$ for some $0\leq k \leq s,$ with $s=nN-2n-1$. This is the same as having a relation
$$R_m: a_0f_0+a_1f_1+....a_nf_n=0$$
where  $a_j \in S$ are homogeneous of degree
$m=s-k$ and $R_m$ is not a consequence of the relations
\begin{equation} 
\label{tij} 
T_{ij}:f_jf_i-f_if_j=0.
\end{equation} 
Since $p_j$ is not a singularity for $D$ for $j>q$, it follows that $f_0(p_j)\ne 0$ in this range.
Hence, for $j>q$, the relation $R_m$ implies that the germ of function induced by $a_0$ at $p_j$ (dividing by some homogeneous polynomial $b_j$ of degree $m$ such that $b_j(p_j)\ne 0$) belongs to the ideal in $\OO_{p_j}$ spanned by the local equations of the divisors $D_1,...,D_n$.  

We apply now the Cayley-Bacharach Theorem as stated in \cite{EGH}, Theorem CB7.

Let $\Gamma$ be the $0$-dimensional subscheme of $\PP^n$ defined by the partial derivatives $f_1,...,f_n$.
Let $\Gamma'$ and $\Gamma''$ be subscheme of $\Gamma$, residual to one another in $\Gamma$, and such that
the support of $\Gamma'$ is the set $\NN'=\{p_{q+1},...,p_r\}$ and the support of $\Gamma''$ is the set $\NN$.
Intuitively, $\Gamma'$ is the 'restriction' of the scheme $\Gamma$ to $\NN'$ and $\Gamma''$ is the 'restriction' of the scheme $\Gamma$ to $\NN$. In particular, the scheme $\Gamma''$ is reduced.

Note that the above discussion implies that the dimension of the family of hypersurfaces $a_0$ of degree $m=s-k$ containing $\Gamma'$ (modulo those containing all of $\Gamma$) is exactly the dimension of $H^{n}(K^*(f))_{nN-n-1-k}$.

On the other hand, for $s$ as above and $0\leq k\leq s$, the Cayley-Bacharach Theorem says that this dimension is equal to the defect $\defect S_k(\NN)$, thus proving the first claim in Theorem \ref{onenode}.

Next we have $\dim H^{n}(K^*(f))_j=\dim H^{n+1}(K^*(f))_{j+N}- \dim H^{n+1}(K^*(f_s))_{j+N} 
=\dim M(f)_{j+N-n-1}-\dim M(f_s)_{j+N-n-1}.$
Moreover, $j \geq n(N-1)$ is equivalent to $j+N-n-1>(n+1)(N-2)$ and hence $\dim M(f)_{j+N-n-1}=\tau(D)=|\NN|$ and $\dim M(f_s)_{j+N-n-1}=0$, thus proving the second claim in Theorem \ref{onenode}.

\endproof

\begin{ex}
\label{exonenode}

We use the notation from Theorem \ref{onenode} and set $T=(n+1)(N-2)$.

\medskip

\noindent (i) If $|\NN|=1$, then $\defect S_k(\NN)=0$ for $k \geq 0$ and therefore we have $ct(D)=st(D)=T$.

\medskip

\noindent (ii) If $|\NN|=2$, then $\defect S_0(\NN)=1$ and $\defect S_k(\NN)=0$ for $k \geq 1$. It follows that $ct(D)+1=st(D)=T$.

\medskip

\noindent (iii) If $|\NN|=3$, and the three nodes are not collinear, then  $\defect S_0(\NN)=2$ and $\defect S_k(\NN)=0$ for $k \geq 1$. It follows that
$ct(D)+1=st(D)=T$
unless $n=2$ when $ct(D)=st(D)=T-1$.

For three collinear points, $\defect S_0(\NN)=2$, $\defect S_1(\NN)=1$ and $\defect S_k(\NN)=0$ for $k \geq 2$. It follows that
$ct(D)+2=st(D)=T$ and 
$\dim M(f)_{T-1}=n+2.$

\medskip

To have some explicit examples of these two distinct situations, consider the following two curves of degree $N=4$:
$$C:f=(x^3+y^3+z^3)x=0$$
and 
$$C':f'=x^2y^2+y^2z^2+x^2z^2-2xyz(x+y+z)-(2xy+3yz+4xz)^2=0.$$
Then the curve $C$ has $3$ collinear nodes and the corresponding Hilbert-Poincar\'e series is
$$HP(M(f))(t)=1+3t+6t^2+7t^3+6t^4+4t^5+3t^6+...$$
with $st(C)=6$. In fact the coefficients of $t^k$ for $0 \leq k \leq 2N-4=4$ are determined by Theorem
\ref{thm3} and the remaining terms are determined by Theorem \ref{onenode}.

In the same way one may obtain
$$HP(M(f'))(t)=1+3t+6t^2+7t^3+6t^4+3t^5+...$$
with $st(C')=5$, using the fact that $C'$ has  $3$ non-collinear nodes located at $(1:0:0)$, $(0:1:0)$ and $(0:0:1)$.

\medskip

\noindent (iv) Here is one example of a sextic curve with $6$ nodes. Consider the curve
$$C:f=x^2(x+z)^2(x-z)^2-y^2(y-z)^2(y^2+2z^2)=0.$$
Then the curve $C$ has $6$ nodes, three of them on the line $y=0$ (namely $(0:0:1)$, $(1:0:1)$ and $(-1:0:1)$)
and the other three on the line $y-z=0$ (namely $(0:1:1)$, $(1:1:1)$ and $(-1:1:1)$).
The corresponding Hilbert-Poincar\'e series is
$$HP(M(f))(t)=1+3t+6t^2+10t^3+15t^4+18t^5+19t^6+18t^7+15t^8+10t^9+7t^{10}+6t^{11}+...$$
with $st(C)=11$.
This result follows exactly by the same argument as above, using in addition the equalities
$\defect S_0(\NN)=5$, $\defect S_1(\NN)=3$, $\defect S_2(\NN)=1$, and $\defect S_k(\NN)=0$ for $k \geq 3$.

\end{ex}

\begin{rk}
\label{RK}
Let $C:f=0$ be a degree $N$ curve in $\PP^2$ such that any singular point of $C$ which is not a node is a unibranch singularity, and let $\NN$ denote the set of its nodes. Then
$\defect S_k(\NN)=0$ for $k>N-3$ and $\defect S_{N-3}(\NN)=r-1$, where $r$ is the number of irreducible components of $C$.

\medskip

This can be derived as follows.  Let $I'$ be the ideal of functions in $S$ vanishing at the points in $\NN$. Then it is shown in Proposition 3.6 in \cite{RK} that there is a minimal resolution
$$ 0 \to \oplus_{i=1,t} S(-b_i) \to \oplus_{j=1,t+1} S(-a_j) \to S \to S/I' \to 0$$
such that $0<a_j <N$ for any $j$, $0<b_i \leq N$ for all $i$ and 
$$ |\{i~~:~~b_i=N\}|=r-1.$$
In fact, Proposition 3.6 in \cite{RK} is stated only for curves with nodes and ordinary cusps, but the only point in the proof where one uses the ordinary cusps is to derive the equality (10), which may also be obtained in our slightly more general setting from the diagram (3.14) on page 201 in \cite{D1}.

\medskip

The above resolution implies that the Hilbert-Poincar\'e series of $S/I'$ is given by the following
equality
$$ HP(S/I')(t)=\frac{1 - \sum_j t^{a_j}+\sum_i t^{b_i}}{(1-t)^3}.$$
Since $\NN$ is a finite set of points, it follows that this series can be rewritten as
$$ HP(S/I')(t)=\frac{Q(t)}{1-t}$$
where $Q(t)$ is a polynomial in $t$ of degree at most $N-2$, the coefficient $c_{N-2}$ of $t^{N-2}$ being exactly
$r-1$.
It follows that $\dim (S/I')_k=|\NN|$ for $k \geq N-2$ and 
$$\dim (S/I')_{N-3}=|\NN|-c_{N-2}= |\NN|-r+1,$$
which proves our claim since one has $\defect S_k(\NN)=|\NN|-\dim (S/I')_k$ for any $k$.

Alternatively, one may end the proof using the formula for the defect or superabundance $\defect S_k(\NN)$
as the difference between the Hilbert polynomial and the Hilbert function given in \cite{RK},
just before the statement of Lemma 3.4.

\medskip

Note that the other main results of our paper {\it do not} extend to this more general setting.
For instance, the curve $C$ constructed in Example \ref{lowdegrel} (ii) for $p=2$, $q=3$, $N=5$ has as singularities two unibranch singularities located at $(1:0:0)$ and $(0:1:0)$ and has a relation of degree 1, i.e. $H^2(K^*(f))_3 \ne 0$, and hence Theorem \ref{thm5} and its consequence Theorem \ref{thm3} fail in this case. 

Moreover, Example \ref{FP2} shows that the spectral sequences considered in the second section in the presence of even ordinary cusps may have a different behaviour than in the case of nodes. Indeed, the differential $d_1: E_1^{1,0} \to E_1^{2,0}$ is trivial for a nodal curve and it is non trivial in Example \ref{FP2}.

The same example shows that Theorem \ref{onenode} fails also in this more general setting, since
$\dim H^2(K^*(f))_4=1$ and $\defect S_1(\NN)=|\NN|-\dim (S/I')_1=0-0=0$ since $\NN=\emptyset.$ 

The resolutions constructed in \cite{DSt} for the Jacobian ideals of Chebyshev curves show that there are no similar results to Proposition 3.6 in \cite{RK} for such Jacobian ideals, not even for nodal curves.

\end{rk}

\begin{rk}
\label{3folds}

For a nodal $3$-fold $D:f=0$ in $\PP^4$ of degree $N$, the fact that $D$ is factorial (i.e. the quotient $S/(f)$ is a unique factorization domain) can be expressed as a vanishing property of a certain defect, namely $\defect S_{2N-5}(\NN)=0$, see Remark 1.2 in Cheltsov's paper \cite{Che}. It follows that Theorem 1.4 in \cite{Che} can be restated as saying that $\defect S_{2N-5}(\NN)=0$ when $|\NN|<(N-1)^2$, which in turn maybe restated in view of Theorem  \ref{onenode} as saying that the corresponding space of syzygies $R_{2N-4}$ is trivial in such a case. On the other hand, Theorem 2.1 (ii) in \cite{DSt3} implies that $R_m=0$ for $m<2N-4$ and {\it any } nodal $3$-fold $D$ in $\PP^4$ of degree $N$.

\end{rk}

\section{The spectral sequence in the case of a nodal surface} \label{sec5}

Let $S:f=0$ be a nodal surface in $\PP^3$ of degree $N$. Then $S$ is a $\Q$-homology manifold
satisfying $b_0(S)=b_4(S)=1$, $b_1(S)=b_3(S)=0$ and the middle Betti number $b_2(S)$ is computable, e.g. using
the formula $b_2(S)=b_2(S_s)-n(S)$, where $S_s$ is a smooth surface in $\PP^3$ of degree $N$, the corresponding second Betti number is given by
$$b_2(S_s)=\frac{(N-1)^4-1}{N}+2$$
and $n(S)=\tau (S)$ is the number of nodes, i.e. the cardinal of the set $\NN$ of nodes of $S$. It follows that
the complement $U$ has at most two non-zero cohomology groups.
The first of them, $H^0(U)$ is $1$-dimensional and of Hodge type $(0,0)$, so nothing interesting here.
The second one, $H^3(U)$, is dual to $H^3_c(U)(-3)$ and $H^3_c(U)$ is isomorphic to $\coker (H^2(\PP^3) \to H^2(S))$, the morphism being induced by the inclusion $i:S \to \PP^3$. It follows that the MHS on $H^3(U)$ is pure of weight $4$
with  
$$h^{4,0}(H^3(U))=h^{0,4}(H^3(U))=0,$$ 
$$h^{3,1}(H^3(U))=h^{1,3}(H^3(U))=h^{2,0}(S)=h^{2,0}(S_s)=p_g(S_s)$$
where the geometric genus of $S_s$ is given by
$$p_g(S_s)={N-1 \choose 3}$$
and
$$h^{2,2}(H^3(U))=h^{1,1}(S)-1=h^{1,1}(S_s)-n(S)-1.$$
In particular we have $P^1H^3(U)=F^1H^3(U)=H^3(U)$ as in Remark \ref{rkA}.

\medskip

Now we look at the nonzero terms in the $E_1$-term of the spectral sequence $E_r^{p,q}(f)$.
They are sitting on two lines, given by $L: p+q=3$ and $L': p+q=2$.

\medskip

We look first at the terms on the line $L$. The term $E_1^{3,0}(f)=H^{4}(K^*)_{N}$ is isomorphic as a $\C$-vector space to $M(f)_{N-4}$, hence has dimension $p_g=p_g(S_s)$.
Hence, the corresponding limit term $E^{3,0}_{\infty}(U)=P^3H^{2}(U)$ has the  dimension at most $p_g$. On the other hand, the above formulas for $h^{p,q}(H^3(U)$ imply that $\dim F^3H^{2}(U)=p_g$. We conclude by Theorem \ref{thmDD} that $F^3H^{3}(U)=P^3H^{3}(U)$ in this case.

\medskip

The term $E_1^{2,1}(f)=H^{4}(K^*)_{2N}$ is isomorphic  to $M(f)_{2N-4}$. Theorem \ref{onenode} 
and Griffiths' results in the smooth case in \cite{Gr}  imply that
\begin{equation} 
\label{dim2.2}
\dim M(f)_{2N-4}=M(f_s)_{2N-4}+ \defect S_{2N-4}(\NN)=h^{1,1}(S_s)-1+ \defect S_{2N-4}(\NN).
\end{equation}
The term $E_1^{1,2}(f)=H^{4}(K^*)_{3N}$ is isomorphic to $M(f)_{3N-4}$, and hence
\begin{equation} 
\label{dim2.3}
\dim M(f)_{3N-4}=M(f_s)_{3N-4}+ \defect S_{N-4}(\NN)= p_g+ \defect S_{N-4}(\NN) .
\end{equation}

\medskip

The term $E_1^{3-q,q}(f)=H^{4}(K^*)_{(q+1)N}$ for $q\geq 3$ is isomorphic to $M(f)_{(q+1)N-4}$, which has dimension $n(S)$.
Moreover, the corresponding limit terms $E^{3-q,q}_{\infty}(U)=Gr_P^{3-q}(H^{2}(U))$ clearly vanish for $q\geq 3$.

\medskip

We look now at the terms on the line $L'$. It follows that the term $E_1^{2,0}(f)=H^{3}(K^*)_{N}$ has dimension
$ \defect S_{2N-4}(\NN)$,
by Theorem \ref{onenode}.
On the other hand $E_{\infty}^{2,0}(f)=0$ which implies in view of the equality $E_{\infty}^{3,0}(f)=
E_{1}^{3,0}(f)$ established above that in fact $ \defect S_{2N-4}(\NN)=0$.

The dimension of the term $E_1^{1,1}(f)=H^{3}(K^*)_{2N}$ is equal to $ \defect S_{N-4}(\NN)$,
again by Theorem \ref{onenode}. And again $E_{\infty}^{1,1}(f)=0$ implies in view of the equality $E_{\infty}^{3,0}(f)=
E_{1}^{3,0}(f)$ established above that the differential $d_1: E_1^{1,1}(f) \to E_1^{2,1}(f)$ is injective.

\medskip

The dimension of the other terms $E_1^{2-q,q}(f)=H^{3}(K^*)_{(q+1)N}$ for $q\geq 2$ is equal to  $n(C)$ and
the corresponding  differentials $d_1: E_1^{2-q,q}(f) \to E_1^{3-q,q}(f)$ are injective
by Theorem  \ref{thm2}, (ii).

In this way we have proved the following.

\begin{thm}
\label{thm7} Let $S:f=0$ be a nodal surface in $\PP^3$ of degree $N$ and let $\NN$ denote the set of its nodes. Then the following hold.

\medskip

\noindent (i) The $E_1$-spectral sequences $E^{p,q}_r(U)$, $E^{p,q}_r(\tilde A)$  and $E^{p,q}_r(f)$
degenerate at the $E_2$-term, i.e. $E_2=E_{\infty}$ for any of these $E_1$-spectral sequences.

\medskip

\noindent (ii) The subspace $P^3H^3(U)=F^3H^3(U)$ of $H^3(U)$ has dimension $p_g={N-1 \choose 3}$.

\medskip

\noindent (iii) $\dim Gr_P^2(H^3(U))=h^{1,1}(S_s)-1-\defect S_{N-4}(\NN)$ and 
$\dim Gr_F^2(H^3(U))=h^{1,1}(S_s)-1-n(S).$ In particular, $P^2H^3(U)=F^2H^3(U)$ if and only if the nodal surface $S$ is smooth or $N<4$.

\medskip

\noindent (iv) $ \defect S_{2N-4}(\NN)=0$.

\end{thm}

\begin{rk}
\label{rkFP}

\noindent (i)
Let $I$ be the homogeneous ideal in $S$ of polynomials vanishing on the set of nodes $\NN$.
The above formulas imply that
$$\dim Gr_F^2(H^3(U))=\dim (I/J_f)_{2N-4}$$
which is a special case of Theorem 2.2 in \cite{DSW}.

\noindent (ii) The ideal $I$ defined in $(i)$ occurs also in the following formula, again a consequence of
Theorem \ref{thm7}.
$$\dim P^2H^3(U)-\dim F^2H^3(U)=\dim (S/I)_{N-4}.$$
When the number of nodes is large, this difference can also be very large. For instance, if $S$ is a Chebyshev surface in $\PP^3$ whose affine equation is 
$$T_N(x)+T_N(y)+T_N(z)+1=0,$$ 
with $T_N(t)$ the degree $N$ Chebyshev polynomial in $\C[t]$, then $I_{N-4}=0$, see Proposition 3.1 in \cite{DSt3}.
It follows that in this case
$$\dim P^2H^3(U)-\dim F^2H^3(U)=\dim S_{N-4}= {N-1 \choose 3}.$$

\end{rk}

\end{document}